\documentclass[smallextended,envcountsame]{svjour3}
\smartqed  
\usepackage{mathptmx}      
\usepackage[dvips]{graphicx,psfrag}
\usepackage{amsmath,amssymb}
\usepackage{color}
\newcommand{\Ref}[1]{{\rm (\ref{#1})}}
\newcommand{\R}{\mathbb R}
\newcommand{\bfx}{\mathbf x}

\newcommand{\It}{{I}_\tau}
\newcommand{\Itn}{{I}_{\tau_n}}
\newcommand{\Ik}{{I}_{K}}
\newcommand{\PP}{\mathcal{P}}
\newcommand{\dd}{\mathrm{d}}
\journalname{Japan J.\ Indust.\ Appl.\ Math.}

\begin{document}

\title{On the Circumradius Condition \\
   for Piecewise Linear Triangular Elements}
\author{Kenta Kobayashi, Takuya Tsuchiya}

\titlerunning{The circumradius condition}

\author{Kenta Kobayashi \and Takuya Tsuchiya}


\institute{Kenta Kobayashi \at
           Graduate School of Commerce and Management, \\
         Hitotsubashi University, Japan \\
           \email{kenta.k@r.hit-u.ac.jp} \and
           Takuya Tsuchiya \at
           Graduate School of Science and Engineering, \\
           Ehime Univesity, Japan \\
           \email{tsuchiya@math.sci.ehime-u.ac.jp}
}

\date{Received: date / Accepted: date}

\maketitle

\begin{abstract}
We discuss the error analysis of linear interpolation
on triangular elements.  We claim that the circumradius condition
is more essential than the well-known maximum angle condition
for convergence of the finite element method, 
especially for the linear Lagrange finite element.
Numerical experiments show that this condition is the best possible.
 We also point out that the circumradius condition is
closely related to the definition of surface area. 
\keywords{linear interpolation \and the circumradius condition
\and the finite element method 
\and Schwarz's example \and the definition of surface area}
\subclass{65D05, 65N30, 26B15}
\end{abstract}

\section{Introduction}
In numerical analysis, linear interpolation on triangular elements is
one of the more fundamental conceptions. Specifically, as meshes become
finer, it is an important tool in understanding why and how finite
element approximations converge to an exact solution.

Let $\Omega \subset \R^2$ be a bounded polygonal domain. Suppose that
we would like to solve the Poisson equation of finding
$u \in H_0^1(\Omega)$ for a given $f\in L^2(\Omega)$ such that
\begin{equation}
    - \Delta u = f \quad \text{ in } \quad \Omega.
   \label{poisson}
\end{equation}
With a triangulation $\tau$ of $\Omega$, we define
the FEM solution $u_h$ by
\begin{equation}
   \int_{\Omega} \nabla u_h\cdot \nabla v_h \dd x = 
   \int_\Omega f v_h \dd x, \qquad \forall v_h \in S_\tau,
  \label{fem}
\end{equation}
where $\PP_1$ is the set of all polynomials whose degree is at most $1$
and
\[
    S_\tau := \left\{v_h \in H_0^1(\Omega)
   \bigm| v_h|_K \in \PP_1, \forall K\in \tau\right\}.
\]
Let $K \subset \R^2$ be a triangle with apices $\bfx_i$, $i = 1, 2, 3$.
We shall always consider $K$ to be a closed set in $\R^2$.
For a continuous function $f \in C^0(K)$,
 the linear interpolation $\Ik f \in \PP_1$ is defined by
\[
   (\Ik f)(\bfx_i) = f(\bfx_i), \qquad i = 1, 2, 3.
\]
If $f \in C^0(\overline{\Omega})$, the linear interpolation
$\It f$ is defined by $(\It f)|_K = \Ik f$.
C\'ea's lemma claims that the error $|u - u_h|_{1,2,\Omega}$
is estimated as
\begin{align*}
  |u - u_h|_{1,2,\Omega} \le \inf_{v_h \in S_\tau}
   |u - v_h|_{1,2,\Omega} \le |u - \It u|_{1,2,\Omega}
   = \left(\sum_{K\in\tau} |u - \Ik u|_{1,2,K}^2 \right)^{1/2}.
\end{align*}
Therefore, the interpolation error $|u - \Ik u|_{1,2,K}$
provides an \textit{upper bound} of $|u - u_h|_{1,2,\Omega}$.

It has been known that we need to impose a geometric condition
to $K$ to obtain an error estimation of $|u - \Ik u|_{1,2,K}$.
We mention the following well-known results.  Let $h_K$ be the
diameter of $K$ and $\rho_K$ be the maximum radius of the inscribed
circle in $K$.
\begin{itemize}
\item \textbf{The minimum angle condition}, Zl\'amal \cite{Z} (1968). \\
{\it Let $\theta_0$, $0 < \theta_0 < \pi/3$, be a constant. If
any angle $\theta$ of $K$ satisfies $\theta \ge \theta_0$
and $h_K \le 1$, then there exists a constant $C = C(\theta_0)$
independent of $h_K$ such that}
\[
   \|v - \Ik v\|_{1,2,K} \le C h_K |v|_{2,2,K}, \qquad
   \forall v \in H^2(K).
\]
 \item \textbf{The regularity (inscribed ball) condition},
 see, for example, Ciarlet \cite{Ciar}. \\
{\it Let $\sigma > 0$ be a constant. If $h_K/\rho_K \le \sigma$ and
$h_K \le 1$, then there exists a constant $C = C(\sigma)$ independent
of $h_K$ such that}
\[
   \|v - \Ik v\|_{1,2,K} \le C h_K |v|_{2,2,K}, \qquad
   \forall v \in H^2(K).
\]
\item \textbf{The maximum angle condition},
 Babu\v{s}ka-Aziz \cite{BA}, Jamet \cite{J} (1976). \\
{\it Let $\theta_1$, $2\pi/3 \le \theta_1 < \pi$, be a constant.  If
any angle $\theta$ of $K$ satisfies $\theta \le \theta_1$ and
$h_K \le 1$, then there exists a constant $C = C(\theta_1)$ independent
of $h_K$ such that}
\[
   \|v - \Ik v\|_{1,2,K} \le C h_K |v|_{2,2,K}, \qquad
   \forall v \in H^2(K).
\]
\end{itemize}
It is easy to show that the minimum angle condition is equivalent
to the regularity condition \cite[Exercise~3.1.3, p130]{Ciar}.
Since its discovery, the maximum angle condition was believed
to be the most essential condition for convergence of solutions of
the finite element method.

 However, Hannukainen-Korotov-K\v{r}\'{i}\v{z}ek pointed out that
\textit{``the maximum angle condition is not necessary for
convergence of the finite element method''} by showing simple
numerical examples \cite{HKK}. 
We double checked the first numerical experiment in \cite{HKK}
with slightly different triangulations and obtained the same result
for the error associated with the finite element approximations.
Therefore, the question arises:
\textit{``What is the essential condition to impose on triangulation
for convergence of the finite element method?''}.
One of the aims here is to give a partial answer to this
question.  Suppose that a sequence $\{\tau_n\}_{n=1}^\infty$ of
triangulations of $\Omega$ is given.
Let $R_K$ be the circumradius of a triangle $K$ and
$R_{\tau_n} := \max_{K\in\tau_n} R_K$.  We claim that the condition
\begin{equation}
   \lim_{n\to \infty} R_{\tau_n} = 0
  \label{cr}
\end{equation}
is \textit{more essential} than the maximum angle condition.
The condition \Ref{cr} is called the \textbf{circumradius condition}.

We moreover point out that the circumradius condition  is closely
related to the definition of surface area.  In the 19th century, people
believed that surface area could be defined  as the limit of the area of
inscribed polygonal surfaces.
In the 1880s, Schwarz and Peano independently presented their famous
example that refutes this expectation. See 
\cite{Cesa}, \cite{Rado1}, \cite{Rado2}. 
We shall observe in Section~3 that, in Schwarz's example,
the limit of the inscribed polygonal surfaces is equal to the area of
the cylinder if and only if the circumradius of triangles converges to
$0$.

We shall also show that the graph of $f \in W^{2,1}(\Omega)$ has
finite area $A_L(f)$.  Moreover, the areas of its inscribed polygonal
surfaces converge to $A_L(f)$ if the sequence of triangulations
satisfies the circumradius condition.  See Theorem~\ref{thm1}
in Section~3.

\vspace{0.3cm}
Let us summarize the notation and terminology to be used.
The Lebesgue and Sobolev spaces on a domain $\Omega \subset \R^2$
are denoted by $L^p(\Omega)$ and $W^{m,p}(\Omega)$, $m=1,2$, 
$1 \le p \le \infty$.  As usual, $W^{m,2}(\Omega)$ is denoted by
$H^m(\Omega)$.  The norms and seminorms of $L^p(\Omega)$ and
$W^{m,p}(\Omega)$ are denoted by $\|\cdot\|_{m,p,\Omega}$ and
$|\cdot|_{m,p,\Omega}$, $m = 0, 1, 2$, $1 \le p \le \infty$.
For a polygonal domain $\Omega \subset \R^2$, a \textit{triangulation}
$\tau$ is a set of triangles which satisfies the following
properties: (recall that each $K$ is a closed set.)
\begin{itemize}
 \item[(i)] $\displaystyle \bigcup_{K\in\tau} K = \overline{\Omega}$,
  and   $\mathrm{int}K \cap \mathrm{int}K' = \emptyset$ for
  any $K$, $K'\in\tau$ with $K\neq K'$.
 \item[(ii)] If $K \cap K' \neq \emptyset$ for $K$, $K'\in\tau$,
  $K \cap K'$ is either their apices or their edges.
\end{itemize}
For a triangulation $\tau$, we define
$|\tau|:=\max_{K\in \tau}\mathrm{diam}K$.

Let $\Omega \subset \R^2$ be a bounded Lipschitz domain.
By Sobolev's imbedding theorem, we have the continuous inclusion 
$W^{2,p}(\Omega) \subset C^0(\overline{\Omega})$ for any $p \in [1,\infty]$. 
Note that for $p=1$ Morry's inequality is not applicable
and the inclusion $W^{2,1}(\Omega) \subset C^0(\overline{\Omega})$ is not so
obvious.  For a proof of the critical 
imbedding, see \cite[Theorem~4.12]{AdamsFournier} and
\cite[Lemma~4.3.4]{BrennerScott}.

\section{Kobayashi's formula, the circumradius condition, and
Schwarz's example}
Recently, we made progress on the error analysis of linear
interpolation on triangular elements.
Liu-Kikuchi presented an explicit form of the constant $C$
in the maximum angle condition \cite{LK}.
Being inspired by Liu-Kikuchi's result, Kobayashi,
one of the authors, obtained the following remarkable result with
the assistance of numerical validated computation \cite{K1}.
\begin{theorem}[Kobayashi's formula] \label{thm2}
Let $A$, $B$, $C$ be the lengths of the three edges of $K$
and $S$ be the area of $K$.  Define the constant $C(K)$ by
\[
   C(K) := \sqrt{\frac{A^2B^2C^2}{16S^2} - 
                    \frac{A^2 + B^2 + C^2}{30} - \frac{S^2}5
       \left(\frac1{A^2} + \frac1{B^2} + \frac1{C^2}\right)},
\]
then the following estimate holds:
\[
   |v - \Ik v|_{1,2,K} \le C(K) |v|_{2,2,K}, \qquad
   \forall v \in H^2(K).
\]
\end{theorem}
Let $R_K$ be the radius of the circumcircle of $K$.  From the formula
$R_K = ABC/4S$, we realize $C(K) < R_K$ and obtain a corollary
of Kobayashi's formula.
\begin{corollary}\label{cor1}
For any triangle $K \subset \R^2$, the following estimate holds:
\[
   |v - \Ik v|_{1,2,K} \le R_K |v|_{2,2,K}, \qquad
   \forall v \in H^2(K).
\]
\end{corollary}

Let $\theta_K \ge \pi/3$ be the maximum angle of $K$.
By the law of sines, we have $h_k = 2R_K\sin \theta_K$.
Therefore, if there is a constant $\theta_1$,
$2\pi/3 \le \theta_1 < \pi$ such that $\theta_K \le \theta_1$, then
$h_K \ge (2 \sin\theta_1) R_K$ and $\lim_{h_K \to 0} R_K = 0$.
This means that, under the assumption $h_K \to 0$, (i)
\textit{the maximum angle condition implies the circumradius condition}.

\begin{figure}[thb]
\begin{center}
  \psfrag{h}[][]{$h$}
  \psfrag{k}[][]{$h^\alpha$}
  \includegraphics[width=5cm]{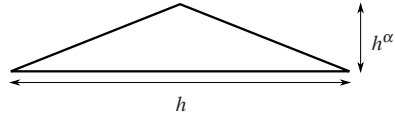}
 \caption{An example of triangles which violates the maximum
angle condition but $R_K \to 0$ as $h\to 0$.}
\label{fig1}
\end{center}
\end{figure}

Let $\Omega \subset \R^2$ be a bounded domain and an arbitrary
$v \in H^2(\Omega)$ is taken and fixed. 
Consider now the isosceles triangle $K \subset \Omega$ depicted in
Figure~\ref{fig1}. If $ 0 < h < 1$ and $\alpha > 1$, then $h^\alpha < h$
and the circumradius of $K$ is $R_K = h^\alpha/2 + h^{2-\alpha}/8$.
Hence, Kobayashi's formula and its corollary yield that, if $1 < \alpha < 2$,
$|v - \Ik v|_{1,2,K} \le R_K |v|_{2,2,K}$ and $R_K \to 0$
as $h\to 0$, whereas the maximum angle of $K$ approaches $\pi$.
This means that, when $h_K \to 0$,
(ii) \textit{the circumradius condition does not necessarily imply the
maximum angle condition}.

Gathering from (i) and (ii), we infer that the circumradius of a
triangle is a more important indicator than its minimum and maximum
angles.

Without the assistance of numerical validated computation,
the authors then proved for arbitrary $p \in [1,\infty]$
the following theorem.

\begin{theorem}[The circumradius condition \cite{KobaTsuchi}]
\label{thm3}
For an arbitrary triangle $K$ with $R_K \le 1$, there
exists a constant $C_p$ independent of $K$
such that the following estimate holds: 
\begin{equation}
   \|v - \Ik v\|_{1,p,K} \le C_p R_K |v|_{2,p,K}, \qquad
   \forall v \in W^{2,p}(K), \quad 1 \le p \le \infty.
     \label{circumradius-inequality}
\end{equation}
\end{theorem}
For the case $p=2$, the estimate \Ref{circumradius-inequality} was shown by
Rand in his Ph.D.\ dissertation \cite[Theorem~7.10]{Rand} but it was not
published in a research paper.

Combining C\'ea's lemma and Corollary~\ref{cor1} or Theorem~\ref{thm3},
we immediately obtain the following estimation.
\begin{theorem} \label{thm5}
Let $u$ be the exact solution
of \Ref{poisson} and $u_h$ be the FEM solution
of \Ref{fem}.  Suppose that $u\in H^2(\Omega)$. Then we have,
for $R_\tau \le 1$,
\begin{equation}
  \|u - u_h\|_{1,2,\Omega} \le C R_\tau |u|_{2,2,\Omega}, \qquad
  R_\tau := \max_{K\in\tau}R_K,
  \label{fem-error}
\end{equation}
where the positive constant $C$ depends only on $C_2$ and $\Omega$.
\end{theorem}
Note that it follows from Corollary~1 that $C_2 = 1$.
However, proving this without using
validated numerical computation is not easy.

\begin{figure}[b]
\begin{center}
   \psfrag{A}[][]{$2\pi r$}
   \psfrag{B}[][]{$H$}
  \includegraphics[width=4.5cm]{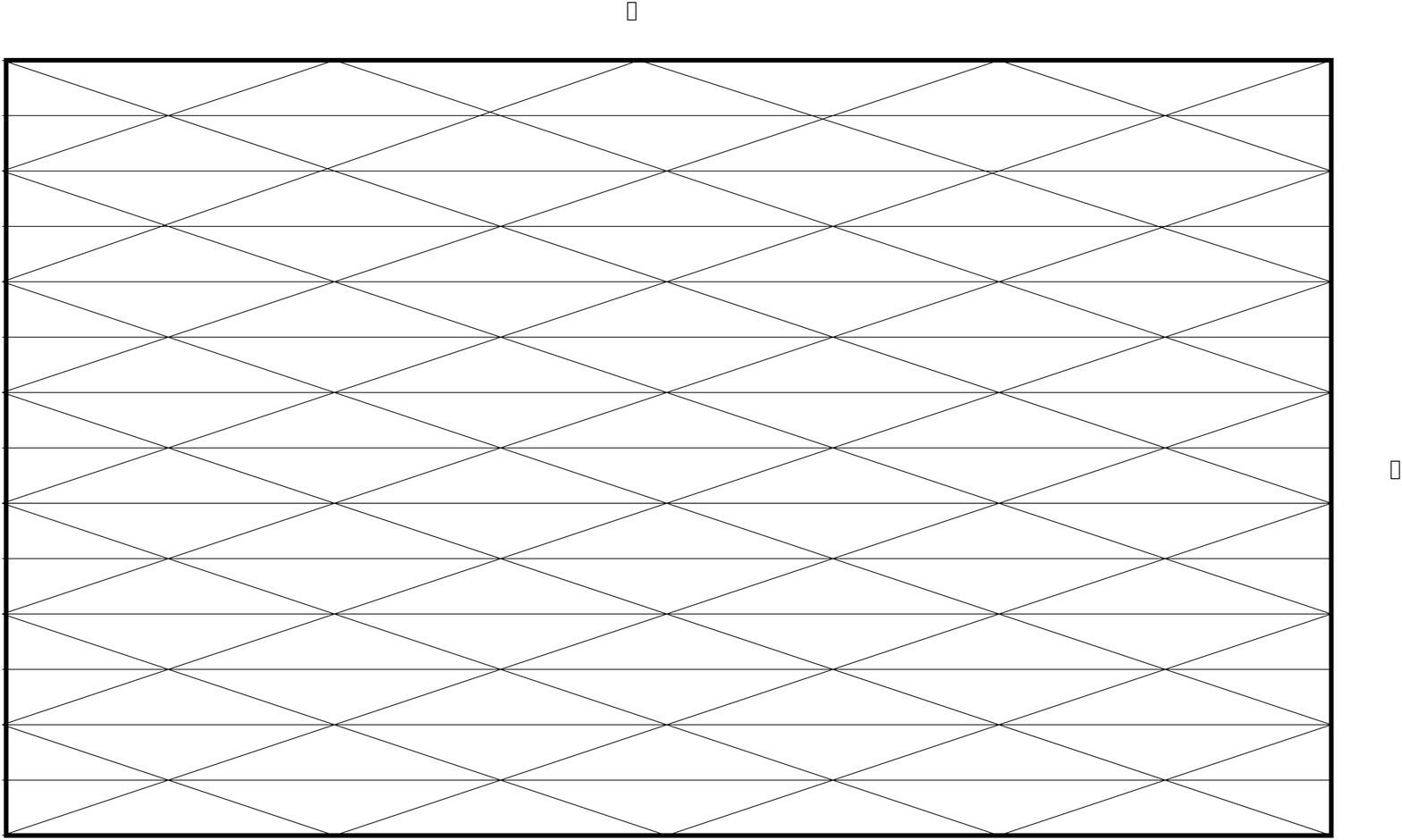}
  \hspace{0.5cm}
  \includegraphics[width=3cm,angle=90]{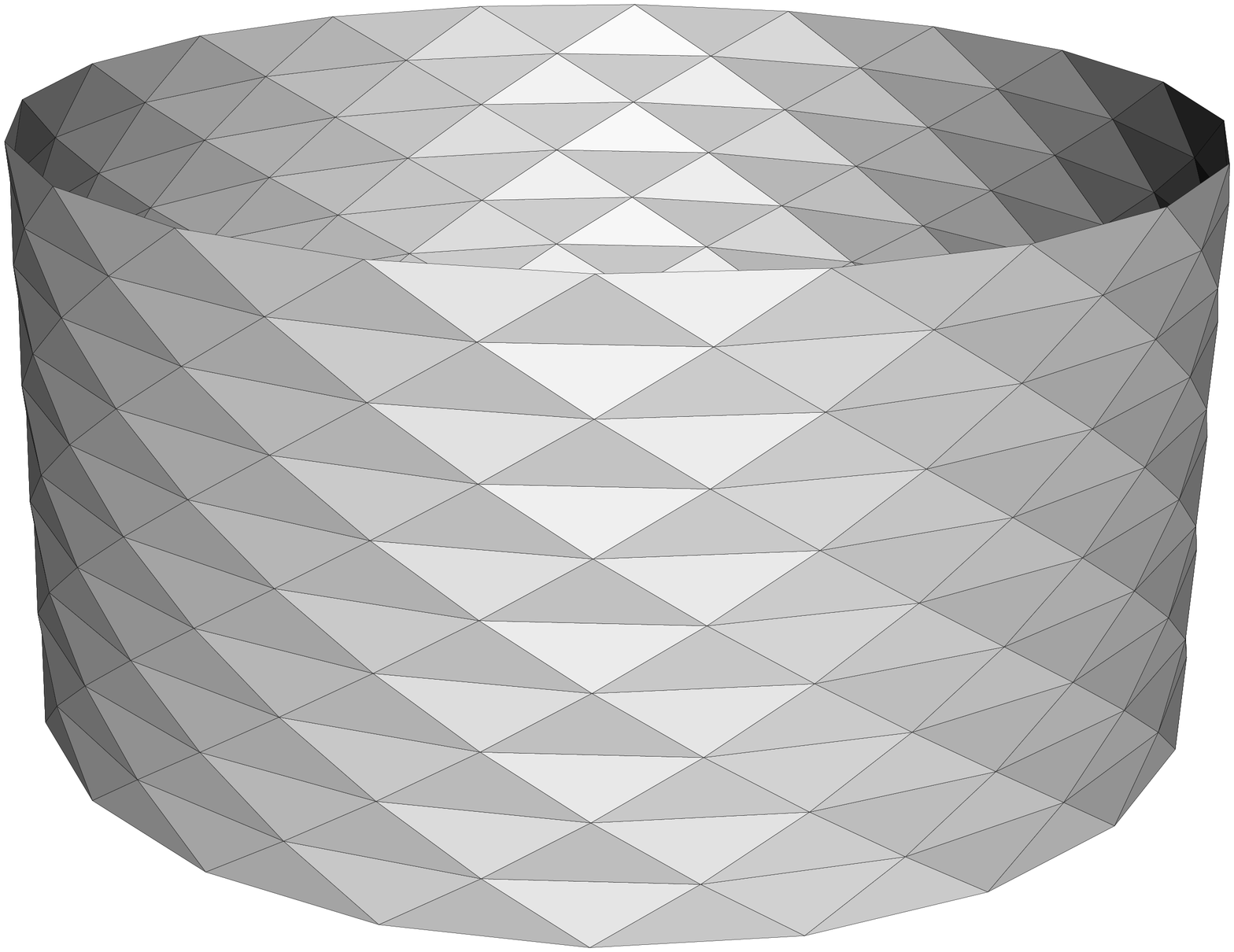} 
 \caption{Schwarz's example.}
 \label{fig2}
\end{center} 
\end{figure}

The isosceles triangle in Figure~1 reminded the authors of Schwarz's
example.  As is well understood, the length of a curve is defined as the
limit of the length of the inscribed polygonal edges.  Hence, one might
think that the area of a surface could be defined in a similar manner.
Actually, mathematicians in the 19th century believed that the area of
surface is the limit of the areas of inscribed polygonal surfaces.

In the 1880s, Schwarz and Peano independently showed, however, that this
definition does not work 
\cite{Cesa}, \cite{Rado1}, \cite{Rado2}, \cite{Zames}.
Let $\Omega$ be a rectangle of height $H$ and width $2\pi r$.
Let $m$, $n$ be positive integers.  Suppose that this rectangle
is divided into $m$ equal strips, each of height $H/m$.
Each strip is then divided into isosceles triangles whose
base length is $2\pi r/n$, as depicted in Figure~\ref{fig2}.
Then, the piecewise linear map $\varphi_\tau:\Omega \to \R^3$
is defined by ``rolling up this rectangle'' so that
all vertexes are on the cylinder of height $H$ and radius $r$.
Then, the cylinder is approximated by the inscribed polygonal surface
which consists of $2mn$ congruent isosceles triangles.
Because the height of each triangle is
$\sqrt{(H/m)^2 + r^2(1 - \cos(\pi/n))^2}$
and the base length is $2r\sin(\pi/n)$, 
the area $A_E$ of the inscribed polygonal surface is
\footnote{The subscript `$E$' of $A_E$ stands for `Elementary'.}
\begin{align*}
  A_E & = 2 m n r \sin \frac{\pi}{n} 
  \sqrt{\left(\frac{H}{m}\right)^2 + r^2 
        \left(1 - \cos \frac{\pi}{n}\right)^2} \\
   & = 2 \pi r \frac{\sin \frac{\pi}{n}}{\frac{\pi}{n}}
  \sqrt{H^2 + \frac{\pi^4r^2}{4} \left(\frac{m}{n^2}\right)^2
   \left(\frac{\sin \frac{\pi}{2n}}{\frac{\pi}{2n}}\right)^4 }.
\end{align*}
If $m,n \to \infty$, we observe
\begin{align*}
& \lim_{m,n \to \infty} A_E = 2\pi r \sqrt{H^2 + \frac{\pi^4 r^2}4 
    \lim_{m,n \to \infty} \frac{m}{n^2}},
\end{align*}
in particular, 
\begin {align*}
  \lim_{m,n \to \infty} A_E = 2\pi r H  \;\; \text{ if and only if } 
  \lim_{m,n \to \infty} \frac{m}{n^2} = 0.
\end{align*}

As we are now aware that the circumradius is an important factor,
we compute the circumradius $R$ of the isosceles triangle in Schwarz's
example. By a straightforward computation, we find that
\begin{align*}
  R = \frac{\frac{H^2}m + \pi^2 r^2 \frac{m}{n^2}
  \left(\frac{\sin \frac{\pi}{2n}}{\frac{\pi}{2n}}\right)^2
   }
  {2 \sqrt{H^2 + \frac{\pi^4r^2}{4}
   \left(\frac{m}{n^2}\right)^2
  \left(\frac{\sin \frac{\pi}{2n}}{\frac{\pi}{2n}}\right)^4}}
\end{align*}
and immediately realize that
\begin{align}
 \lim_{m,n \to \infty} A_E = 2\pi r H \Longleftrightarrow
  \lim_{m,n \to \infty} \frac{m}{n^2} = 0 \Longleftrightarrow
  \lim_{m,n \to \infty} R = 0.
   \label{equiv}
\end{align}
This fact strongly suggests that the circumradius of triangles in
a triangulation is essential for error estimations of linear
interpolations.

\begin{figure}[thb]
\begin{center}
      \includegraphics[width=5.5cm]{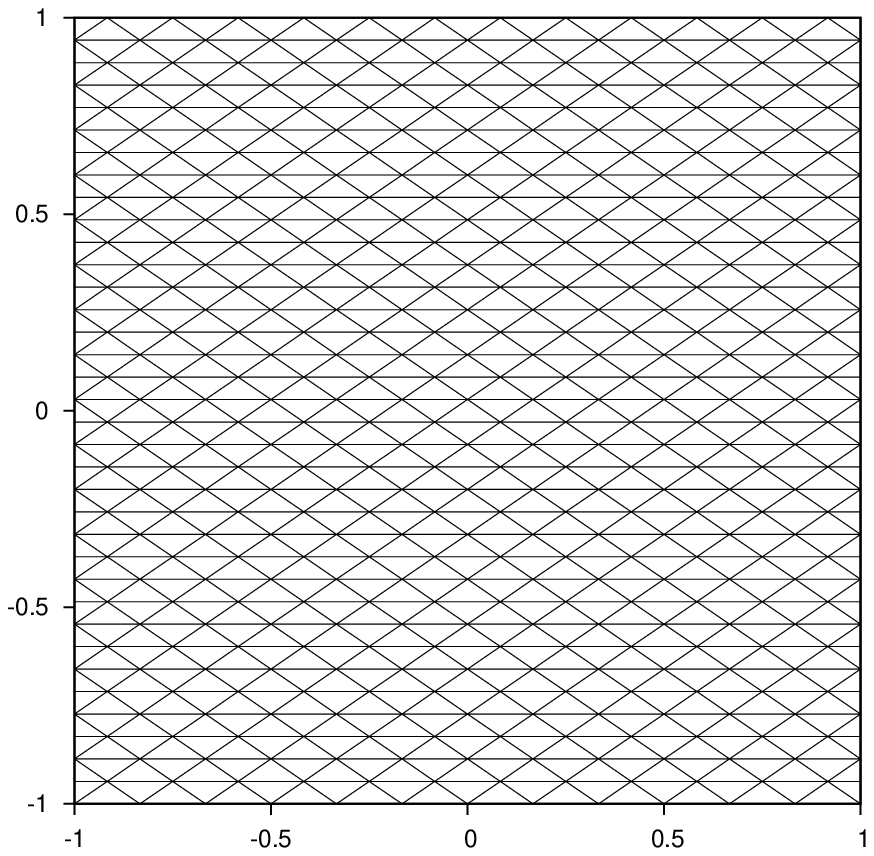}
          \includegraphics[width=5.5cm]{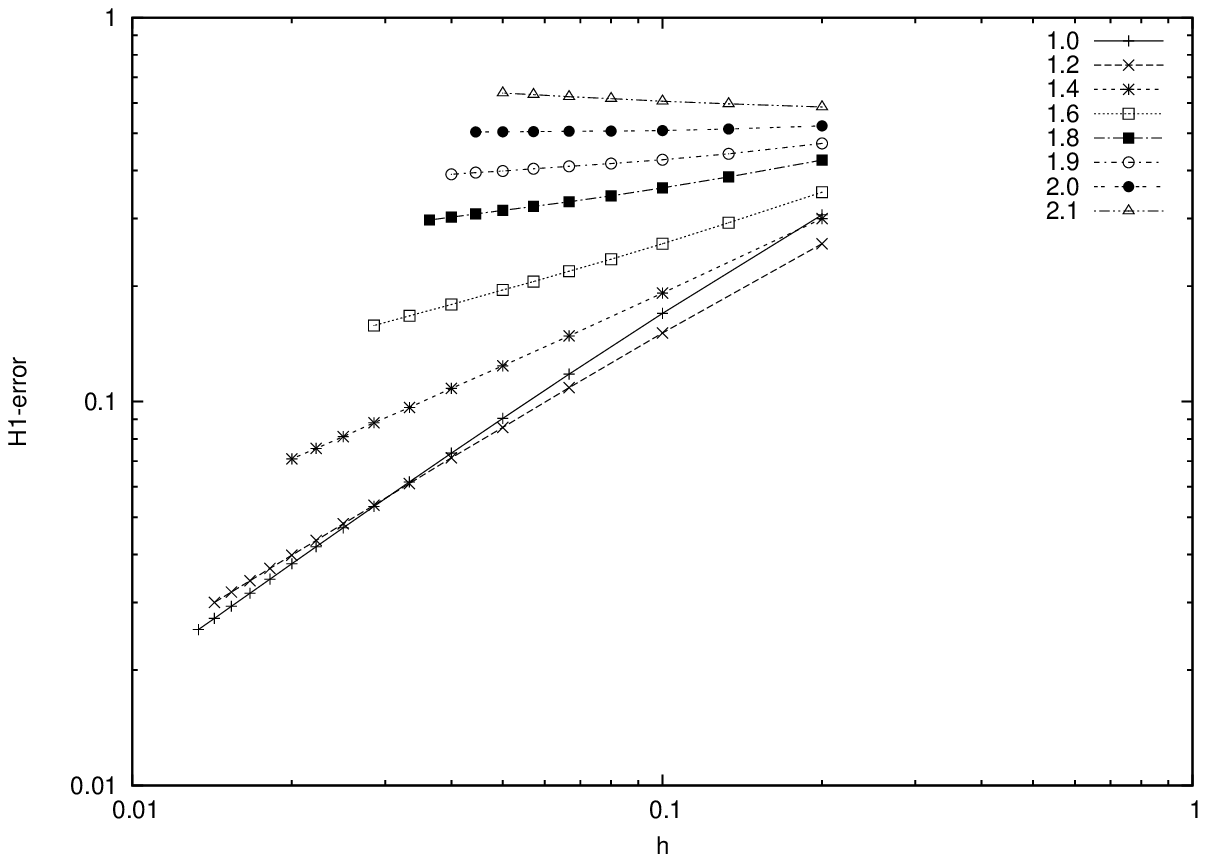}
\end{center}
 \caption{The triangulation of $\Omega$ with $N=12$ and $\alpha=1.6$
  and the errors for FEM solutions in $H^1$-norm.
  The horizontal axis represents the maximum diameter of the
 triangles and the vertical axis represents $H^1$-norm of the errors
 of the FEM solutions. The number next to the symbol
indicates the value of $\alpha$.}
\label{fig3}
\end{figure}

With \Ref{equiv} in mind, we perform a numerical experiment similar
to the one in \cite{HKK}.  Let $\Omega := (-1,1)\times(-1,1)$, 
$f(x,y) := a^2/(a^2 - x^2)^{3/2}$, and
$g(x,y) := (a^2 - x^2)^{1/2}$ with $a:=1.1$.
Then we consider the following
Poisson equation:  Find $u \in H^1(\Omega)$ such that
\begin{equation}
   - \Delta u = f \quad \text{ in } \Omega, \qquad
   u = g \quad \text{ on } \partial\Omega.
  \label{test_problem}
\end{equation}
The exact solution of \Ref{test_problem} is $u(x,y) = g(x,y)$
and its graph is a part of the cylinder.
For a given positive integer $N$ and $\alpha > 1$, we consider the
isosceles triangle with base length $h:=2/N$ and height
$2/\lfloor 2/h^\alpha\rfloor \approx h^\alpha$, as depicted
in Figure~\ref{fig1}.  For comparison, we also consider the isosceles
triangle with base length $h$ and height $h/2$ for 
$\alpha=1$.  We triangulate $\Omega$
with this triangle, as shown in Figure~\ref{fig3}.
The behavior of the error is given in Figure~3.  The horizontal axis
represents the mesh size measured by the maximum diameter of triangles
in the meshes and the vertical axis represents the error associated with
FEM solutions in $H^1$-norm.  The graph clearly shows that
the convergence rates worsen as $\alpha$ approaches $2.0$.
For $\alpha=2.1$, the FEM solutions even diverge.  We replot the
same data in Figure~4, in which the horizontal axis represents the
maximum of the circumradius of triangles in the meshes.  Figure~4 shows
convergence rates are almost the same in all cases if we measure these
with the circumradius.
\begin{figure}[thb]
\begin{center}
    \includegraphics[width=5.5cm]{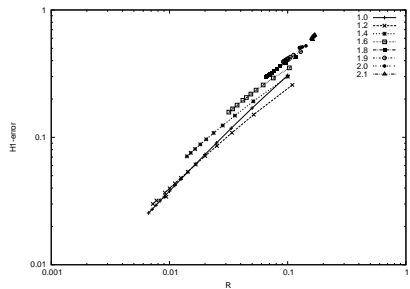}
\end{center}
 \caption{Replotted data: the errors in $H^1$-norm of FEM
solutions measured using the circumradius. The horizontal axis represents
the maximum circumradius of the triangles.}
\label{fig4}
\end{figure}

From the results of the numerical experiments, we draw the following
conclusions: Suppose that we consider the Poisson equation \Ref{test_problem}.
\begin{itemize} 
\item In our example, although the triangulation does not satisfy the
      maximum angle condition, the FEM solutions converge to the exact
      solution and the error behaves exactly as the estimation
      \Ref{fem-error} predicts. If the triangulation does not satisfy
      the circumradius condition, the FEM solutions diverge even if
      meshes become finer with respect to the maximum diameter of the
      triangles.  From this observation, we infer that, for convergence
      of the FEM solutions, the circumradius solution is
      \textsl{more essential} than the maximum angle condition and
      is the \textsl{best possible} as a geometric condition for
      triangulation.
      \footnote{
      By the statement (i) given after
      Corollary~\ref{cor1}, we realize that if the circumradius
      condition does not hold then the maximum angle condition does not
      hold either.}
 \item 
      The numerical experiments in \cite{HKK} show that,
      in certain combinations of an exact solution and triangulation,
      FEM solutions can converge to an
       exact solution, although triangulation does not satisfy
       the maximum angle condition.  We notice that their
       triangulations do not satisfy the circumradius condition either.
      Hence, the circumradius condition is not necessary
      for convergence of the finite element method.
\end{itemize}
These conclusions answer, partially but not completely, the question
which Hannukai-nen-Korotov-K\v{r}\'{i}\v{z}ek posed.
We infer from the numerical experiments that matching between
exact solutions and geometry of triangulation seems important.
Further and deeper understanding of how FEM solutions converge
to an exact solution is strongly desired.

\section{The circumradius condition and the definition of surface area}
At the present time, the most general definition of surface area is 
that of Lebesgue.  Let $\Omega:=(a,b)\times (c,d) \subset \R^2$ be a
rectangle and $\tau_n$ be a sequence of triangulation of $\Omega$ such
that $\lim_{n \to \infty}|\tau_n| = 0$.  Let
$f \in C^0(\overline{\Omega})$ be a given continuous function.
Let $f_n\in S_{\tau_n}$ be such that $\{f_n\}_{n=1}^\infty$ converges
uniformly to $f$ on $\overline{\Omega}$.  Note that the graph of
$z = f_n(x,y)$ is a set of triangles and its area is defined as a sum of
these trianglar areas.  We denote this area by $A_E(f_n)$ and have
\[
   A_E(f_n) = \int_{\Omega} \sqrt{1 + |\nabla f_n|^2} \dd x.
\]
Let $\Phi_f$ be the set of all such sequences
 $\{(f_n,\tau_n)\}_{n=1}^\infty$.
Then the area $A_L(f) = A_L(f;\Omega)$ of the graph $z=f(x,y)$
is defined by
\[
   A_L(f) = A_L(f;\Omega):= \inf_{\{(f_n,\tau_n)\} \in \Phi_f} 
  \liminf_{n\to \infty}A_E(f_n).
\]
This $A_L(f)$ is called the \textbf{surface area of} $z = f(x,y)$ 
\textbf{in the Lebesgue sense}.
For a fixed $f$, $A_L(f;\Omega)$ is additive and continuous
with respect to the rectangular domain $\Omega$.
Tonelli then presented the following theorem.

For a continuous function $f\in C^0(\overline{\Omega})$, we define
$W_1(x)$, $W_2(y)$ by
\begin{align*}
 W_1(x) := \sup_{\tau(y)} \sum_{i} |f(x,y_{i-1}) - f(x,y_i)|,
   \quad x \in (a,b),\\
 W_2(y) := \sup_{\tau(x)} \sum_{j} |f(x_{j-1},y) - f(x_j,y)|,
   \quad y \in (c,d),
\end{align*}
where $\tau(y)$, $\tau(x)$ are subdivisions
 $c = y_0 < y_1 < \cdots < y_N = d$ and
 $a = x_0 < x_1 < \cdots < x_M = b$, respectively and `$\sup$' are
taken for all such subdivisions.
Then, a function $f$ has 
\textbf{bounded variation in the Tonelli sense} if
\[
   \int_a^b W_1(x) \dd x + \int_c^d W_2(y) \dd y < \infty.
\]
Also, a function $f$ is called 
\textbf{absolutely continuous in the Tonelli sense} if, 
for almost all $y \in (c,d)$ and $x \in (a,b)$,
the functions $g(x) := f(x,y)$ and $h(y) := f(x,y)$
are absolutely continuous on $(a,b)$ and $(c,d)$, respectively.

\begin{theorem}[Tonelli] \label{Tonelli}
For a continuous function $f \in C(\overline{\Omega})$ defined on a
rectangular domain $\Omega$, its graph $z = f(x,y)$ has
finite area $A_L(f) < \infty$ if and only if $f$ has bounded variation
in the Tonelli sense. If this is the case, we have
\begin{equation}
   A_L(f) \ge \int_\Omega \sqrt{1 + f_x^2 + f_y^2}\, \dd \bfx.
  \label{tonelli}
\end{equation}
In the above inequality, the equality holds if and only if
$f$ is absolutely continuous in the Tonelli sense.
\end{theorem}

For a proof of this theorem, see \cite[Chapter V, pp.163--185]{Saks}.
It follows from Tonelli's theorem that if $f \in W^{1,\infty}(\Omega)$
then the area $A_L(f)$ is finite and the equality holds in
\eqref{tonelli}.  In the following theorem we consider the
case $f \in W^{2,1}(\Omega)$.

\vspace{-0.1cm}
\begin{theorem} \label{thm1}
Let $\Omega \subset \R^2$ be a rectangular domain.
If $f \in W^{2,1}(\Omega)$, then its graph has finite area,
that is, $A_L(f) < \infty$, and 
the equality holds in \eqref{tonelli}.
Moreover, if a sequence $\{\tau_n\}_{n=1}^\infty$ of triangulations 
of $\Omega$ satisfies the circumradius condition, then we have
\[
   \lim_{n\to \infty} \int_\Omega \sqrt{1 + |\nabla \Itn f|^2}
\dd x = A_L(f) =
  \int_\Omega \sqrt{1 + |\nabla f|^2}\, \dd \bfx.
\]
\end{theorem}
\noindent
\textbf{Proof.} 
At first, we notice $f$ is of bounded variation and absolutely
continuous in the Tonelli sense.
Let $\omega:= \{x\}\times(c,d)$ for $x \in (a,b)$.
We consider the trace operator
$\gamma:W^{2,1}(\Omega) \to W^{1,1}(\omega)$ defined by
$(\gamma f)(y):=f(x,y)$.  Then, $\gamma$ is a bounded linear
operator and it is easy to see that
\begin{align*}
  \sum_{\tau(y)} |(\gamma f)(y_{i-1}) - (\gamma f)(y_i)|
  \le \int_c^d |(\gamma f)'(y)| \dd y, \quad
  W_1(x) \le \int_c^d |f_y(x,y)| \dd y.
\end{align*}
Similarly, we obtain $W_2(y) \le \int_a^b |f_x(x,y)| \dd x$ and
Fubini's theorem implies that $f$ has bounded variation in the
Tonelli sense. Hence, Theorem~\ref{Tonelli} yields $A_L(f) < \infty$ and
\eqref{tonelli} holds.
We show that $f$ is absolutely continuous in the Tonelli sense
in exactly the same manner.
Therefore, the equality holds in \eqref{tonelli}.

For the piecewise linear interpolation $\Itn f$, we have
\begin{align*}
  A_E(\Itn f) = A_L(\Itn f) = 
  \int_\Omega\sqrt{1 + (\Itn f)_x^2 + (\Itn f)_y^2} \, \dd \bfx.
\end{align*}
Hence, $|A_L(f) - A_E(\Itn f)|$ is estimated as
{\allowdisplaybreaks
\begin{align*}
  |A_L(f) - &  A_E(\Itn f)|  \le \int_\Omega
  \left| \sqrt{1 + f_x^2 + f_y^2} -
   \sqrt{1 + (\Itn f)_{x}^2 + (\Itn f)_{y}^2} \right| \dd \bfx \\
  & \le \int_\Omega
  \frac{\left|\left(f_x + (\Itn f)_{x}\right)
   \left(f_x - (\Itn f)_{x}\right) + 
   \left(f_y + (\Itn f)_{y}\right)
   \left(f_y - (\Itn f)_{y}\right)\right|}
  {\sqrt{1 + f_x^2 + f_y^2} +
   \sqrt{1 + (\Itn f)_{x}^2 + (\Itn f)_{y}^2}}
   \dd \bfx \\
 & \le |f - \Itn f|_{1,1,\Omega} \\
& \le C_1 R_{\tau_n} |f|_{2,1,\Omega}
  \to 0 \qquad \text{ as } R_{\tau_n} \to 0,
\end{align*}
}
because
\begin{align*}
   & \frac{|f_x + (\Itn f)_{x}|}{\sqrt{1 + f_x^2 + f_y^2} +
   \sqrt{1 + (\Itn f)_{x}^2 + (\Itn f)_{y}^2}} \le 1, \\
   & \frac{|f_y + (\Itn f)_{y}|}{\sqrt{1 + f_x^2 + f_y^2} +
   \sqrt{1 + (\Itn f)_{x}^2 + (\Itn f)_{y}^2}} \le 1.
\end{align*}
Thus, Theorem~\ref{thm1} is proved.  $\square$

\vspace{0.3cm}
Note that, from Schwarz's example,  Theorem~\ref{thm1} is the
\textit{best possible} with respect to the geometric condition
for triangulation.  At this point, one might be tempted
to define the surface area using the circumradius condition in the
following way:

\begin{definition}
Let $\Omega \subset \R^2$ be a bounded polygonal domain.
Suppose that a sequence $\{\tau_n\}$ of triangulation of $\Omega$
satisfies the circumradius condition.  Then, for a continuous function
$f \in C^0(\overline{\Omega})$, the area $A_{CR}(f)$ of the surface
$z = f(x,y)$ is defined by
\[
   A_{CR}(f) := \lim_{n \to \infty} A_E(\Itn f).
\]
\end{definition}

Theorem~\ref{thm1} claims that, for $f \in W^{2,1}(\Omega)$,
$A_{CR}(f)$ is well-defined and $A_{CR}(f)=A_L(f) < \infty$.
The example given by Besicovitch shows that $A_{CR}(f)$ is not
well-defined in $C^0(\overline{\Omega})$ in general
\cite{Besicovitch}.  That is, there exists
$f \in C^0(\overline{\Omega})$ and two triangulation sequences 
$\{\tau_n\}_{n=1}^\infty$,  $\{\mu_n\}_{n=1}^\infty$ 
of $\Omega$ which satisfy the circumradius
condition such that $A_L(f) < \infty$ and
\[
  \lim_{n \to \infty} A_E(\Itn f) \neq \lim_{n \to \infty}
   A_E({I}_{\mu_n} f).
\]
Therefore, we present the following problem.
Let $X$ be a Banach space such that
$W^{2,1}(\Omega) \subset X \subset L^1(\Omega)$.

\begin{problem}
\begin{itemize} 
 \item[$(1)$] Determine the largest function space $X$ such that
$A_{CR}(f)$ is well-defined for any $f \in C^0(\overline{\Omega}) \cap X$.
\item[$(2)$] With $X$ defined in (1), \textit{prove} or
\textit{disprove} whether $A_{CR}(f) = A_L(f)$ for any
$f \in  C^0(\overline{\Omega}) \cap X$ with $A_L(f) < \infty$.
\end{itemize}
\end{problem}

\section{Concluding remarks ---  \textit{History repeats itself}}
We have shown that the circumradius condition
is more essential than the maximum angle condition for convergence
of FEM solutions.  Also, we have pointed out a close
relationship between the circumradius condition and the definition
of surface area.  In concluding, we draw readers'
attention to the similarity of two histories.
After Schwarz and Peano found their counter example, 
mathematicians naturally tried to find a \textit{proper} definition
of surface area.  The authors are unfamiliar with 
the history behind that quest.  Instead, we suggest
that readers look at \cite[Capter~I]{Rado2} from which
we mention the following remarks.

Let $S \subset \R^3$ be a general parametric surface.  If there
exists a Liptschitz map $\varphi:\Omega \to \R^3$ defined on
a domain $\Omega\subset\R^2$ such that $S = \varphi(\Omega)$,
$S$ is called \textbf{rectifiable}.  
Let $\Omega$ be a rectangle and $\varphi:\Omega \to \R^3$
be a rectifiable surface. Suppose that we have a sequence
$\{\tau_n\}_{n=1}^\infty$ of triangulation of $\Omega$ such
that $|\tau_n| \to 0$ as $n\to\infty$.  Then, the rectifiable
surface $\varphi$ has linear interpolations $I_{\tau_n}\varphi$.
Rademacher showed \cite{Radem1}, \cite{Radem2} that if
$\{\tau_n\}_{n=1}^\infty$ satisfies the minimum angle condition we
have $\lim_{n\to\infty} A_E(I_{\tau_n}\varphi) = A_L(\varphi)$. 
Then, Young showed \cite{Young} that if
$\{\tau_n\}_{n=1}^\infty$ satisfies the maximum angle condition we
have $\lim_{n\to\infty} A_E(I_{\tau_n}\varphi) = A_L(\varphi)$.
See also the comment by Fr\'echet \cite{Frechet} on Young's result.
\footnote{In \cite[p.12]{Rado2}, Rad\'o wrote wrongly that the second
result was by Fr\'echet.}

This means that the minimum and maximum angle conditions were
already found about 50 years before they were rediscovered by
FEM exponents.  This is an interesting example of
the proverb \textit{History repeats itself}.

\begin{acknowledgements}
The first author is supported by Inamori Foundation and
 JSPS Grant-in-Aid for Young Scientists (B) 22740059.
The second author is partially supported by JSPS
Grant-in-Aid for Scientific Research (C) 22440139 and Grant-in-Aid for
Scientific Research (B) 23340023.
The authors thank the anonymous referee for valuable comments
and for drawing the authors' attention to Besicovitch's paper
\cite{Besicovitch}.
\end{acknowledgements}

\end{document}